\documentclass[12pt]{amsart}
\usepackage{amssymb,latexsym,amsmath,eufrak,amsthm,cite,pictex}
\usepackage[active]{srcltx}
\usepackage{hyperref}
\usepackage[inline]{enumitem}
\usepackage[T1]{fontenc}
\usepackage{mathrsfs} 
\usepackage{tikz}
\usetikzlibrary{arrows, calc, patterns}

\usetikzlibrary{arrows.meta}
\usetikzlibrary{positioning}
\usetikzlibrary{bending}
\tikzset{>=triangle 45}
\allowdisplaybreaks[1]

\def\g{\gamma}

\def\s{\sigma}

\def\G{\Gamma}

\def\z{\zeta}

\def\Cs{Ces\'{a}ro}

\newcommand{\C}{{\mathbb C}}

\makeatletter
\newcommand{\inlineitem}[1][]{%
\ifnum\enit@type=\tw@
    {\descriptionlabel{#1}}
  \hspace{\labelsep}
\else
  \ifnum\enit@type=\z@
       \refstepcounter{\@listctr}\fi
    \quad\@itemlabel\hspace{\labelsep}
\fi}
\makeatother

\theoremstyle{plain}
\newtheorem{theorem}{Theorem}

\theoremstyle{definition}

\theoremstyle{remark}

\def\Cs{Ces\`{a}ro}
\title{Remark Concerning  {\Cs} Operator on the Hardy Space $H^p(\C_+)$ in the Upper Half-Plane}
\author{Valentin~V.~Andreev, Miron~B.~Bekker, Joseph~A.~Cima}
\date{}
\begin{document}
\maketitle
\begin{abstract}
 We consider {\Cs} operator on the Hardy space $H^p(\C_+)$ in the upper half-plane for $1<p<\infty$. In \cite{AS} it was proved that for all $1<p<\infty$
 the spectrum of the operator $V=\frac{2(p-1)}{p}C-I$ is located on the unit circle and in \cite{ABC1} the authors of this note showed that for $p=2$ operator $V$ is unitary. In the present note we show that for $1<p<\infty$, $p\ne 2$, the norm of the operator $V$ is strictly greater than one.
\end{abstract}

The {\Cs} operator on the Hardy space $H^p(\C_+)$ in the upper half-plane was introduced in \cite{AS}. It is defined by the expression
\begin{equation}\label{Cesaro-definition}
 (Cf)(z)=\frac{1}{z}\int_0^z f(\z)d\z, \quad f\in H^p(C_+).
\end{equation}
In \cite{AS} was proved that for $1<p\le\infty$ the operator $C$ is bounded on $H^p(\C_+)$ and for $1<p<\infty$ its spectrum $\s(C)$ is given by
\begin{equation}\label{Cesaro-spectrum}
 \s(C)=\{w:|w-\frac{p}{2(p-1)}|=\frac{p}{2(p-1)}\}.
\end{equation}
Authors of \cite{AS} also proved that the point spectrum of the operator $C$ is empty.

Define an operator $V$ on $H^p(C_+)$ by the formula
\begin{equation}\label{V-definition}
 V=\frac{2(p-1)}{p}C-I.
\end{equation}
By the Spectral Mapping theorem the spectrum of $V$ is the unit circle,
\begin{equation}\label{V-spectrum}
 \s(V)=\{w:|w|=1\}.
\end{equation}
In \cite{ABC1} authors of the present note proved that for $p=2$ the operator $V$ is unitary on $H^2(C_+)$. Further properties of the {\Cs} operator on
$H^2(\C_+)$ are investigated in \cite{ABC2}.

Since the spectral radius of $V$ is equal to one for all $p\in(1,\infty)$, we have $\Vert V\Vert_p\ge 1$, where $\Vert V\Vert_p$ means the norm of $V$ considered as an operator on $H^p(\C_+)$.
\begin{theorem}
 Let the operator $V$ on $H^p(\C_+)$ ($1<p<\infty$) be defined by formulas \eqref{Cesaro-definition} and \eqref{V-definition}. If $p\ne 2$ then $\Vert V\Vert_p>1$.
\end{theorem}
{\bf Proof}. To prove the theorem we exibit a vector $f$ which belongs to $H^p(\C_+)$ for all $p\in(1,\infty)$ and show that for all such $p$ $\Vert Vf\Vert_p\ge \Vert f\Vert_p$ with equality only for $p=2$.

Put
\begin{equation*}
 f_t(z)=\frac{-it}{(-it-z)^2},\quad t>0.
\end{equation*}
The function $f_t$ belongs to $H^p$ and
\begin{equation}\label{ft p-norm}
 \Vert f_t\Vert^p_p=\frac{1}{t^{p-1}}\int\limits_{-\infty}^{\infty}\frac{1}{(x^2+1)^p}dx=\frac{1}{t^{p-1}}\frac{\sqrt{\pi}\G(-\dfrac{1}{2}+p)}{\G(p)}.
\end{equation}

Then
\begin{equation*}
 (Cf_t)=\frac{1}{-it-z}
\end{equation*}
and
\begin{equation*}
 (Vf_t)(z)=\frac{2(p-1)}{p}\frac{-i\gamma t-z}{(-it-z)^2},
\end{equation*}
where $\gamma=\g(p)=\dfrac{2-p}{2(p-1)}$. Simple calculations give
\begin{equation}\label{Vft pnorm}
 \Vert Vf_t\Vert^p_p=\frac{1}{t^{p-1}}\left[\frac{2(p-1)}{p}\right]^p\int\limits_{-\infty}^{\infty}\frac{(x^2+\g^2)^{p/2}}{(x^2+1)^p}dx.
\end{equation}
Since $\g(p)=0$ only for $p=2$ it follows that
\begin{gather*}
  \Vert Vf_t\Vert^p_p\ge \frac{1}{t^{p-1}}\left[\frac{2(p-1)}{p}\right]^p2\int\limits_0^{\infty}\frac{x^p}{(x^2+1)^p}dx=\\
  \frac{1}{t^{p-1}}2\sqrt{\pi}\left(\frac{p-1}{p}\right)^{p-1}\frac{\G(\dfrac{p}{2}+\dfrac{1}{2})}{\G(\dfrac{p}{2}+1)}
\end{gather*}
with equality only for $p=2$. Therefore
\begin{equation*}
 \Vert Vf_t\Vert_p^p/\Vert f_t\Vert^p\ge 2\left(\frac{p-1}{p}\right)^{p-1}\frac{\G(\dfrac{p}{2}+\dfrac{1}{2})\G(p)}{\G(\dfrac{p}{2}+1)\G(p-\dfrac{1}{2})}
\end{equation*}
Put
\begin{equation}\label{Phi}
\Phi(p)= \frac{2(p-1)^{p-1}}{p^{p-1}}\frac{\G\left(\dfrac{1}{2}+\dfrac{p}{2}\right)\G(p)}{\G\left(\dfrac{p}{2}+1\right)\G\left(-\dfrac{1}{2}+p\right)}.
\end{equation}
Then
\begin{gather*}
 \Phi(p)=\frac{4(p-1)^{p-1}}{p^p}\frac{\G(\dfrac{p}{2}+\dfrac{1}{2})\G(p)}{\G(\dfrac{p}{2})\G(p-\dfrac{1}{2})}=
 \frac{2}{\sqrt{\pi}}(p-1)^{p-1}\left(\frac{2}{p}\right)^p\frac{\G^2(\frac{p}{2}+\frac{1}{2})}{\G(p-\frac{1}{2})}.
\end{gather*}
Taking logarithmic derivative of both sides one obtains
\begin{equation*}
 \frac{\Phi^{\prime}(p)}{\Phi(p)}=\ln{\left(2(1-\frac{1}{p})\right)}+\frac{\G^{\prime}(\frac{p}{2}+\frac{1}{2})}{\G(\frac{p}{2}+\frac{1}{2})}-\frac{\G^{\prime}(p-\frac{1}{2})}{\G(p-\frac{1}{2})}.
\end{equation*}
Using well known representation of logarithmic derivativa of gamma-function one obtains
\begin{equation}\label{log-der}
  \frac{\Phi^{\prime}(p)}{\Phi(p)}=\ln{\left(2(1-\frac{1}{p})\right)}+(1-\frac{p}{2})\sum\limits_{n=0}^{\infty}\frac{1}{(n+p-1/2)(n+p/2+1/2)}.
\end{equation}
We see, in particular, that $\Phi^{\prime}(2)=0$ and $\Phi^{\prime}(p)\to-\infty$ as $p\to 1^+$. The series in the right side of \eqref{log-der} converges uniformly in $p$ for $p>1$, therefore it can be differentiated term-by-term. We have
\begin{gather*}
 \frac{\Phi^{\prime\prime}(2)}{\Phi(2)}=\frac{d}{dp}\left( \frac{\Phi^{\prime}(p)}{\Phi(p)}\right)(2)=\frac{1}{2}-\frac{1}{2}\sum\limits_{n=0}^{\infty}\frac{1}{(n+3/2)^2}=\\
 \frac{1}{2}\left(1-4\sum\limits_{n=0}^{\infty}\frac{1}{(2n+3)^2}\right)=\frac{1}{2}\left(5-\frac{\pi^2}{2}\right),
\end{gather*}
where we used the well known formula
\begin{equation*}
 \sum\limits_{n=1}^{\infty}\frac{1}{(2n-1)^2}=\frac{\pi^2}{8}.
\end{equation*}
Thus, $\Phi^{\prime\prime}(2)>0$ and $p=2$ is a minimum of the function $\Phi$. From \eqref{Phi} it follows that $\Phi(2)=1$. Further analysis indicates that the function $\Phi$ monitonically increases for $p>2$.
\begin{figure}[!h]
 \begin{center}
  \includegraphics[width=5in]{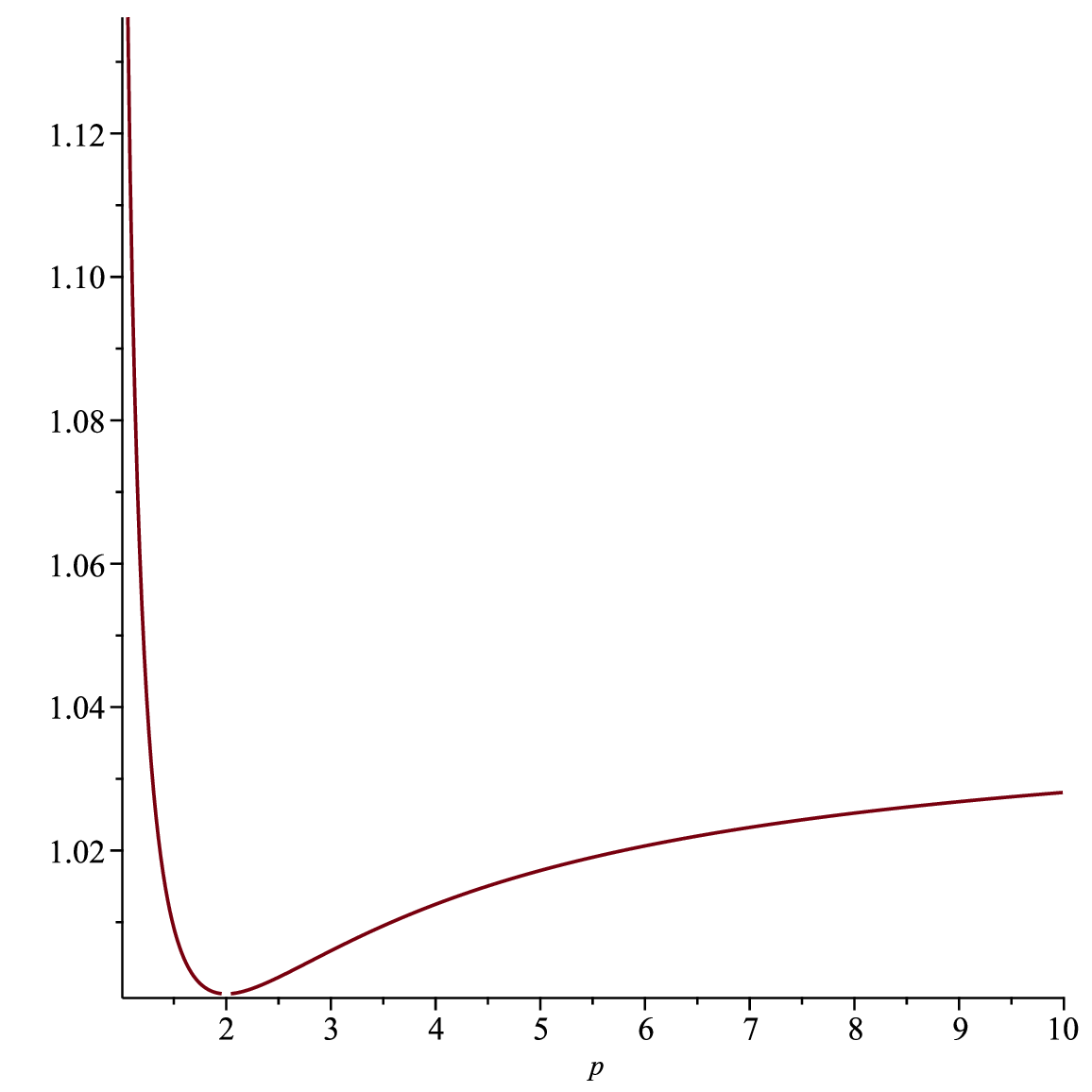}
 \end{center}
\caption{\small $\Phi(p)$ for $1.05\le p\le 10$}
\end{figure}

VALENTIN~V.~ANDREEV Department of Mathematics, Lamar University, Beaumont,TX USA; {\it vvandreev8@gmail.com}\\
MIRON B. BEKKER Department of Mathematics, the University of \\Pittsburgh at Johnstown,
Johnstown, PA, USA, {\it bekker@pitt.edu}\\
JOSEPH~A.~CIMA Department of Mathematics, The University of North Carolina at Chapel Hill, Phillips Hall, Chapel Hill,NC, USA, {\it cima@email.unc.edu}
\end{document}